\documentclass{article}
\usepackage{amsfonts,amsmath,amsthm,amssymb,graphics,epsfig}

\input{tcilatex}
\begin{document}

\title{Singular Points of Real Quintic Curves Via Computer Algebra}
\author{David A. Weinberg and Nicholas J. Willis}
\maketitle


\bigskip

\abstract{There are 42 types of real singular points
for irreducible real quintic curves and 49 types of real singular
points for reducible real quintic curves. The classification of
real singular points for irreducible real quintic curves is
originally due to Golubina and Tai. There are 28 types of singular
points for irreducible complex quintic curves and 33 types of
singular points for reducible complex quintic curves. We derive
the complete classification with proof by using the computer
algebra system Maple.  We clarify that the classification is based
on computing just enough of the Puiseux expansion to separate the
branches. Thus, a major component of the proof consists of a sequence of large symbolic
computations that can be done nicely using Maple.}

\section{Introduction}

The classification of singular points of real quintic curves is originally
due to Golubina and Tai[4]. They found 41 individual types of real singular
points for real irreducible quintic curves. In this paper, we will exhibit,
with proof, the classification of individual types of singular points for
both irreducible and reducible real quintic curves. We found 42 individual
types of singular points for irreducible quintic curves. We think that class
39 of Golubina and Tai should split into two distinct classes based on
whether the two tangent lines are real and distinct or complex conjugate
(these classes are represented by diagrams 3 and 4 below). We exhibit the 49
individual types of singular points for reducible real quintic curves. There
are 28 individual types of singular points for irreducible complex quintic
curves and 33 individual types of singular points for reducible complex
quintic curves. \ Our description of the equivalence relation is new and our
proof is new and gives a very nice illustration of the role that computer
algebra can play in doing proofs. \ Furthermore, our proof is self-contained
and is the most elementary proof possible, which makes the material
accessible to the widest possible audience.

The classification of singular points of complex projective quintic curves
appears in a paper by A. Degtyarev[3]. He not only exhibits the 28
individual types of singular points for complex irreducible quintics, but he
exhibits all 221 sets of singular points, and furthermore proves that the
rigid isotopy type of an irreducible complex quintic is defined uniquely by
its set of singular points.

The general question is how shall we classify singular points of real
quintic curves. For each fixed degree n, we want a finite classification of
singular points for all algebraic curves of degree n. Thus, in general, the
local diffeomorphism type is not the desired criterion of classification.
For example, in the Arnol'd notation, four lines intersecting at the origin
represents an $X_{9}$ singular point, which is really an infinite family of
smoothly inequivalent singularities. \ Notice here that an irreducible real
quintic curve can have an $X_{9}$ singular point. The tradition is to treat
these as one class by fiat. \ In our scheme the $X_{9}$ family will appear
naturally as a single class.

Now let us describe how we will classify the individual types of singular
points that a real quintic curve can have. Given any polynomial equation $%
f(x,y)=0$, it is possible to solve for $y$ in terms of $x$ in the form of
fractional power series, called Puiseux expansions. There is an algorithm
for doing this, and the software Maple computes such Puiseux expansions,
even for curves with literal coefficients. Our classification will be based
on taking just enough of the Puiseux expansions to separate the
\textquotedblleft branches," and noting the exponents at which the
\textquotedblleft branches" separate. In other words, compute the Puiseux
expansions to a power of $x$ such that all expansions are unique. Then we
will associate a tree-type graph, to which we will refer as a "tree diagram"
or "diagram." \ These diagrams will be described in detail below and will
codify how the \textquotedblleft branches" separate and will serve to
classify the type of the singular point. \ At this point let us remark that
the term "branch" already has a traditional meaning in this context. \ We
are really interested in the distinct Puiseux expansion. \ In [9], C. T. C.
Wall has coined the term "pro-branch" for the distinct Puiseux expansions.
It follows from Section 10 of Milnor's book[7] that such a classification
gives a finite number of types for each fixed degree.

In studying a singular point of an algebraic curve, the first thing to look
at is the Newton polygon. (Our Newton polygons will follow the style of
Walker [8].) Corresponding to each segment of the Newton polygon, there is a
quasihomogeneous polynomial [p.195, 2]. If all such quasihomogeneous
polynomials have no multiple factors, then the Newton polygon already tells
us the type of the singularity. (Note that in this case, we know right away
the exponents at which all of the Puiseux expansions separate.) But if there
is a multiple factor, then it is necessary to examine the situation more
closely. For this, we turn to the Puiseux expansions. \ As indicated above,
the relevant definition on which the classification is based is new and
appeals to the Puiseux expansion in an invariant way.

Let us note that we will classify the \textit{real} singular points. (It is
possible for a real quintic curve to have a complex conjugate pair of
singular points. We will avoid this case.) By a simple translation of axes,
we may assume that the singular point is at the origin. We will treat both
irreducible and reducible curves, but note that the notions of irreducible
and reducible are with respect to the complex numbers. Note also that we
will not study reducible curves with multiple components.

The objects being classified are pairs whose first coordinate is a real
quintic curve, specified by a 5th degree polynomial with real coefficients,
considered up to a real nonzero multiplicative constant, and the second
coordinate is a singular point of the curve in the first coordinate. \ Let
the quintic curve be given by $f(x,y)=0$, where \newline
\newline
$%
f(x,y)=a_{00}+a_{10}x+a_{01}y+a_{20}x^{2}+a_{11}xy+a_{02}y^{2}+a_{30}x^{3}+a_{21}x^{2}y+a_{12}xy^{2}+a_{03}y^{3}+a_{40}x^{4}+a_{31}x^{3}y+a_{22}x^{2}y^{2}+a_{13}xy^{3}+a_{04}y^{4}+a_{50}x^{5}+a_{41}x^{4}y+a_{32}x^{3}y^{2}+a_{23}x^{2}y^{3}+a_{14}xy^{4}+a_{05}y^{5}
$.\newline
\newline
Since we may assume that our singular point is at $(0,0)$, we have $a_{00}=0$%
. Since the point is singular, $a_{10}=a_{01}=0$. In this paper we will use
the term "tangent cone" to refer to the terms of lowest degree in $f(x,y)$.
The degree of these terms is called the multiplicity of the point. If the
point is of multiplicity five, then the curve must be reducible since any
homogeneous polynomial of degree $5$ must factor. Thus, for irreducible
curves, we only need to study points of multiplicity four, three, or two.

Let us now explain how all of the cases are enumerated. First we choose the
tangent cone by choosing the tangent lines together with their
multiplicities. \ For the cases of interest, the choice of tangent lines can
be fixed by a linear change of coordinates. \ Moreover, by rotation of axes,
we may assume no tangent line is vertical. For each tangent cone, we
consider all possible Newton polygons. For each Newton polygon, we first
consider the case where none of the quasihomogeneous polynomials
corresponding to the segments of the Newton polygon have a multiple factor.
Then we consider the cases where there is a multiple factor. When there is a
multiple factor, the choice of this factor can be fixed by a linear change
of coordinates. \ In these cases, Maple is used to save many hours of hand
calculation to compute the Puiseux expansion of the corresponding families
of curves. The different types of singular points are then determined by the
vanishing or nonvanishing of certain polynomials in the coefficients of the
families of curves; these polynomials are given to us by the Maple
computations in the form of discriminant-like polynomials in the
coefficients of the Puiseux expansions. The details of this are carried out
in the next section.

Let us now discuss the issue of verifying the existence of irreducible
curves that have a given type of singular point. Observe that for a given
degree, the irreducible curves form a dense open subset in the Zariski
topology on the space of all curves of that degree. When a segment of the
Newton polygon contains a multiple quasihomogeneous factor, we use Maple to
determine the different types of singular points corresponding to that
family, and in this process, a sequence of polynomial conditions on the
coefficients (which turn out to be discriminants) is obtained. With respect
to the Zariski topology, if an irreducible curve is found at any stage of
the sequence, then all prior stages contain irreducible curves. Most of the
time it is obvious that a certain family contains an irreducible
representative. If it is not obvious, then Groebner basis techniques can be
used to show that there is an irreducible representative or that every curve
in a given family is reducible (even when Maple will not show this in
response to the 'factor' command). (For more details on this see the Maple
worksheets on the website of David Weinberg [11].)

The details of the following outline will be carried out in the next
section. \ For irreducible quintic curves, by a linear change of coordinates
, as described above, it suffices to consider the following cases,
indicated, for each multiplicity, by choice of tangent cone, number of
Newton polygons, if greater than one, and choice of multiple
quasihomogeneous factors, if applicable.

\bigskip\underline{Multiplicity 4}

$y^{4}$

$y^{3}(y-x)$

$y^{2}(y-x)^{2}$

$(x^{2}+y^{2})^{2}$

$y^{2}(y-x)(y-2x)$

$y^{2}(x^{2}+y^{2})$

$y(y^{2}-x^{2})(y-gx)$

$y(y-x)(x^{2}+y^{2})$

$(x^{2}+y^{2})(x^{2}+4y^{2})$

Short Maple computations are required only for the 3rd and 4th cases above
because of multiple quasihomogeneous factors in the Newton polygon.

\bigskip\underline{Multiplicity 3}

$y^{3}\qquad \qquad \qquad \qquad $(3 Newton Polygons)

$y^{2}(y-x)\qquad \qquad \ \ \ $(4 Newton Polygons)

$(y-x)(y-2x)(y-3x)$

$(y-x)(x^{2}+y^{2})$

Substantial Maple computation is required only for the family $%
y^{2}(y-x)-x^{5}+2x^{3}y+ax^{4}y+bx^{2}y^{2}+cx^{3}y^{2}+dxy^{3}+ex^{2}y^{3}+fy^{4}+gxy^{4}+hy^{5}
$.

\bigskip\underline{Multiplicity 2}

$y^{2}\qquad \qquad \qquad \qquad $(5 Newton Polygons)

$y^{2}-x^{2}$

$y^{2}+x^{2}$

Substantial Maple computation is required for the family $%
(y+x^{2})^{2}+ax^{5}+bx^{3}y+cxy^{2}+dx^{4}y+ex^{2}y^{2}+fy^{3}+gx^{3}y^{2}+hxy^{3}+jx^{2}y^{3}+ky^{4}+lxy^{4}+my^{5}
$. \ Some spectacular factorizations are performed during the course of that
computation and are indicated in the next section.
\bigskip

For reducible quintic curves, Maple computation is used to examine
the cases where an irreducible conic is tangent to an irreducible cubic. \
(The other cases are enumerated by mathematical common sense.)
\bigskip

Given an algebraic curve with a singular point at the origin, let us now
describe how to associate a tree diagram to this singular point once we have
the Puiseux expansions. Each time at least one \textquotedblleft branch"
separates, record the exponent where that happens. Place all such exponents
in a row at the top. For each exponent in the top row, there corresponds a
column of vertices. Each Puiseux expansion corresponds to exactly one vertex
in that column, and those expansions with the same coefficients up to that
exponent correspond to the same vertex. We start with one vertex on the left
corresponding to the power zero. Line segments are drawn connecting the
vertices from left to right, where each polygonal path from left to right
corresponds to Puiseux expansions having the same set of coefficients up to
a given exponent. The diagram stops at the first exponent where each vertex
in that column corresponds to exactly one Puiseux expansion. Notice that
this tree diagram uniquely specifies the singularity type (up to
permutations of vertices within columns) provided that no tangent line at
the origin is vertical. Braces will join pairs of vertices, within a given
column, corresponding to complex conjugate coefficients. In such a case, the
only real solution of the original equation, satisfying the pair of
expansions indicated by the braces, in a small enough neighborhood of the
origin is (0,0).

In [9] C.T.C. Wall uses the term \textquotedblleft pro-branches" to refer to
the distinct Puiseux expansions belonging to a given singular point, and
then defines a notion of \emph{exponent of contact} between two
pro-branches. It follows from Lemma 4.1.1 on page 68 of [9], that the
diagram we assign to a singular point is invariant under a linear change of
coordinates. \newline
\newline
\underline{Example.} $y^{2}=-x^{3}$.\newline
Notice that $y=\pm \;i\;x^{3/2}$, which can also be written as $y=\pm
(-x)^{3/2}$. For each $x<0$, there are two distinct real solutions for $y$.
Hence, the diagram is (without braces!)
\begin{figure}[h]
\scalebox{.33}{\includegraphics[width=60mm]{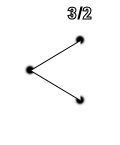}}
\end{figure}
\newline
\newline
\underline{Example.} $x^{2}y+x^{4}+2xy^{2}+y^{3}=0$. \newline
If $B:=x^{2}y+x^{4}+2xy^{2}+y^{3}$, the Maple command puiseux $(B,x=0,y,3)$
tells us that the Puiseux expansions begin as follows:
\begin{equation*}
y=-x+x^{3/2}\qquad (\text{branch}\#1)
\end{equation*}%
\begin{equation*}
y=-x-x^{3/2}\qquad (\text{branch}\#2)
\end{equation*}%
\begin{equation*}
y=-x^{2}\qquad \qquad \;\;(\text{branch}\#3)
\end{equation*}%
In the next section, we will refer to the relevant truncated portion of the
Puiseux expansion as the \emph{Puiseux jet}. Notice that the coefficient of $%
x$ in branch \#1 and branch \#2 is $-1$, while the coefficient of $x$ in
branch \#3 is $0$. So there is a splitting at the first power of $x$, which
is indicated as \newline
\begin{figure}[h]
\scalebox{.33}{\includegraphics{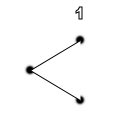}}
\end{figure}
\newline
\newline
\newline
\newline
\newline
Next we must show the splitting of \#1 from \#2. Notice that the power of $x$
at which \#1 and \#2 split is 3/2. Now our diagram looks like\newline
\begin{figure}[h]
\scalebox{.33}{\includegraphics{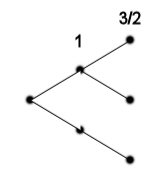}}
\end{figure}
\newline

The diagram is now complete; notice that there are three distinct vertices
in the column labeled $3/2$. \newline

\textbf{Definition of the equivalence relation:} two singular points are
equivalent if they have the same diagram as described above.\newline

To summarize, we have described a precise procedure for assigning a diagram
to a singular point of an algebraic curve and this assignment is invariant
under a linear change of coordinates. \newline
\newline

\underline{Acknowledgments.} The authors wish to thank Tomas Recio
(Universidad de Cantabria, Santander, Spain), Carlos Andradas (Universidad
Complutense de Madrid, Spain), Eugenii Shustin (Tel-Aviv University),
Jeffrey M. Lee (Texas Tech University) and Anatoly Korchagin (Texas Tech
University) for several useful conversations. The authors also wish to thank
Mark van Hoeij (Florida State University) for some singularly valuable Maple
code.


\section{Classification and proof}

\underline{Irreducible curves}\newline
\newline
\underline{Multiplicity 4.}\newline
\newline
Tangent cone: $y^4$\newline
\newline
\begin{tabular}{ll}
\begin{tabular}{l}
Newton polygon: \\
\\
\\
\\
\end{tabular}
& \scalebox{.33}{\includegraphics{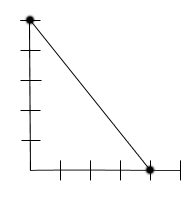}}%
\end{tabular}
\newline
$y^4+ax^5+bx^4y+cx^3y^2+dx^2y^3+exy^4+fy^5=0$, $a\neq0$.\newline
\newline
Puiseux jets:\newline
$y=(-a)^{1/4}x^{5/4}$.\newline
\newline
\begin{tabular}{ll}
\begin{tabular}{l}
1. Diagram: \\
\\
\\
\\
\end{tabular}
& \scalebox{.33}{\includegraphics{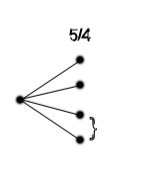}}%
\end{tabular}
\newline
Tangent cone: $y^3(y-x)$\newline
\newline
\begin{tabular}{ll}
\begin{tabular}{l}
Newton polygon: \\
\\
\\
\\
\end{tabular}
& \scalebox{.33}{\includegraphics{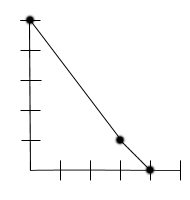}}%
\end{tabular}
\newline
$y^3(y-x)+ax^5+bx^4y+cx^3y^2+dx^2y^3+exy^4+fy^5=0, a\neq0$.\newline
\newline
Puiseux jets:\newline
$y=x$ \newline
$y=(-a)^{1/3}x^{4/3}$.\newline
\newline
\begin{tabular}{ll}
\begin{tabular}{l}
2. Diagram: \\
\\
\\
\\
\end{tabular}
& \scalebox{.33}{\includegraphics{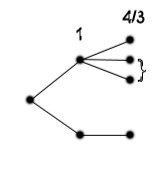}}%
\end{tabular}
\newline
\newline
Tangent cone: $y^2(y-x)^2$\newline
\newline
\begin{tabular}{ll}
\begin{tabular}{l}
Newton polygon: \\
\\
\\
\\
\end{tabular}
& \scalebox{.33}{\includegraphics{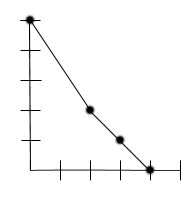}}%
\end{tabular}
\newline
\newline
$A=y^2(y-x)^2+ax^5+bx^4y+cx^3y^2+dx^2y^3+exy^4+fy^5=0, a\neq0$.\newline
\newline
The Puiseux expansion is computed by using the Maple command
puiseux(A,x=0,y,0).\newline
\newline
Puiseux jets: \newline
$y=(-a)^{1/2}x^{3/2}$\newline
$y=x-(a+b+c+d+e+f)^{1/2}x^{3/2}$.\newline
\newline
Condition: $a+b+c+d+e+f\neq0$.\newline
\newline
\begin{tabular}{ll}
\begin{tabular}{l}
3. Diagram: \\
\\
\\
\\
\end{tabular}
& \scalebox{.33}{\includegraphics{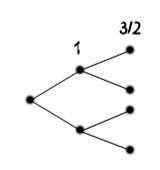}}%
\end{tabular}
\newline
\newline
If $a+b+c+d+e+f=0$, then each curve of the form A is reducible.\newline
\newline
Tangent cone: $(x^2+y^2)^2$.\newline
\newline
\begin{tabular}{ll}
\begin{tabular}{l}
Newton polygon: \\
\\
\\
\\
\end{tabular}
& \scalebox{.33}{\includegraphics{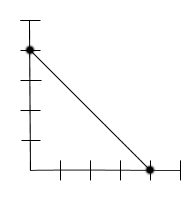}}%
\end{tabular}
\newline
\newline
$B=(x^2+y^2)^2+ax^5+bx^4y+cx^3y^2+dx^2y^3+exy^4+fy^5=0$. \newline
\newline
Conditions: $a-c+e\neq0$ or $b+f-d\neq0$. \newline
\newline
Using Maple, we obtain \newline
\newline
Puiseux jets: \newline
$%
y=RootOf(Z^2+1)x+(a+bRootOf(Z^2+1)+fRootOf(Z^2+1)-c+e-dRootOf(Z^2+1))^{1/2}x^{3/2}
$. \newline
\newline
(Note: 4 expansions here.) \newline
\newline

\begin{tabular}{ll}
\begin{tabular}{l}
4. Diagram: \\
\\
\\
\\
\end{tabular}
& \scalebox{.33}{\includegraphics{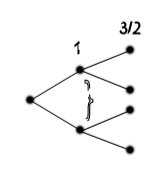}}%
\end{tabular}
\newline
\newline
If $a-c+e=0$ and $b+f-d=0$, then, using Maple, each curve in B is reducible.%
\newline
\newline
Tangent cone: $y^{2}(y-x)(y-2x)$.\newline
\newline
\begin{tabular}{ll}
\begin{tabular}{l}
Newton polygon: \\
\\
\\
\\
\end{tabular}
& \scalebox{.33}{\includegraphics{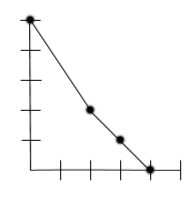}}%
\end{tabular}
\newline
\newline
$y^{2}(y-x)(y-2x)+ax^{5}+bx^{4}y+cx^{3}y^{2}+dx^{2}y^{3}+exy^{4}+fy^{5}=0,a%
\neq 0$.\newline
\newline
Puiseux jets:\newline
$y=x$\newline
$y=2x$\newline
$y=(-\frac{1}{2}a)^{1/2}x^{3/2}$.\newline
\newline
\begin{tabular}{ll}
\begin{tabular}{l}
5. Diagram: \\
\\
\\
\\
\end{tabular}
& \scalebox{.33}{\includegraphics{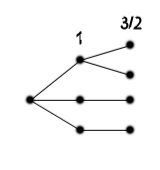}}%
\end{tabular}
\newline
\newline
Tangent cone: $y^{2}(x^{2}+y^{2})$.\newline
\newline
\begin{tabular}{ll}
\begin{tabular}{l}
Newton polygon: \\
\\
\\
\\
\end{tabular}
& \scalebox{.33}{\includegraphics{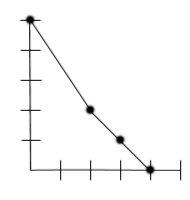}}%
\end{tabular}
\newline
\newline
$%
y^{2}(x^{2}+y^{2})+ax^{5}+bx^{4}y+cx^{3}y^{2}+dx^{2}y^{3}+exy^{4}+fy^{5}=0,a%
\neq 0$.\newline
\newline
Puiseux jets:\newline
$y=\pm ix$\newline
$y=(-a)^{1/2}x^{3/2}$.\newline
\newline
\begin{tabular}{ll}
\begin{tabular}{l}
6. Diagram: \\
\\
\\
\\
\end{tabular}
& \scalebox{.33}{\includegraphics{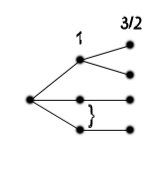}}%
\end{tabular}
\newline
\newline
Tangent cone: $y(y^{2}-x^{2})(y-gx)$.\newline
\newline
\begin{tabular}{ll}
\begin{tabular}{l}
Newton polygon: \\
\\
\\
\\
\end{tabular}
& \scalebox{.33}{\includegraphics{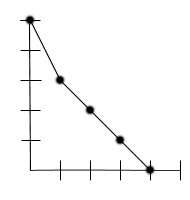}}%
\end{tabular}
\newline
\newline
$%
y(y^{2}-x^{2})(y-gx)+ax^{5}+bx^{4}y+cx^{3}y^{2}+dx^{2}y^{3}+exy^{4}+fy^{5}=0,a\neq 0,g\neq 0,\pm 1
$.\newline
\newline
Puiseux jets: \newline
$y=0x$\newline
$y=x$\newline
$y=-x$\newline
$y=gx$.\newline
\newline
\begin{tabular}{ll}
\begin{tabular}{l}
7. Diagram: \\
\\
\\
\\
\end{tabular}
& \scalebox{.33}{\includegraphics{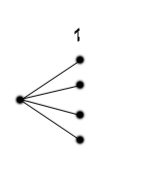}}%
\end{tabular}
\newline
\newline
Tangent cone: $y(y-x)(x^{2}+y^{2})$.\newline
\newline
\begin{tabular}{ll}
\begin{tabular}{l}
Newton polygon: \\
\\
\\
\\
\end{tabular}
& \scalebox{.33}{\includegraphics{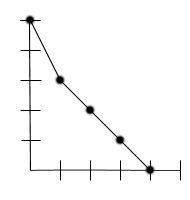}}%
\end{tabular}
\newline
\newline
$%
y(y-x)(x^{2}+y^{2})+ax^{5}+bx^{4}y+cx^{3}y^{2}+dx^{2}y^{3}+exy^{4}+fy^{5}=0,a\neq 0
$.\newline
\newline
Puiseux jets:\newline
$y=0x$\newline
$y=x$\newline
$y=\pm ix$.\newline
\newline
\begin{tabular}{ll}
\begin{tabular}{l}
8. Diagram: \\
\\
\\
\\
\end{tabular}
&
\scalebox{.33}{\includegraphics{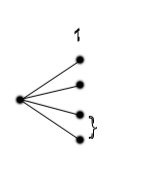}}%
\end{tabular}
\newline
\newline
Tangent cone: $(x^{2}+y^{2})(x^{2}+4y^{2})$.\newline
\newline
\begin{tabular}{ll}
\begin{tabular}{l}
Newton polygon: \\
\\
\\
\\
\end{tabular}
& \scalebox{.33}{\includegraphics{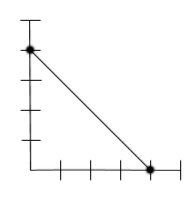}}%
\end{tabular}
\newline
\newline
$%
(x^{2}+y^{2})(x^{2}+4y^{2})+ax^{5}+bx^{4}y+cx^{3}y^{2}+dx^{2}y^{3}+exy^{4}+fy^{5}=0
$.\newline
\newline
Puiseux jets: \newline
$y=\pm ix$\newline
$y=\pm 2ix$.\newline
\newline
\begin{tabular}{ll}
\begin{tabular}{l}
9. Diagram: \\
\\
\\
\\
\end{tabular}
& \scalebox{.33}{\includegraphics{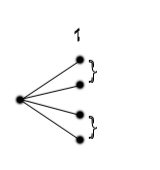}}%
\end{tabular}
\newline
\newline
\underline{Multiplicity 3.}\newline
\newline
Tangent cone: $y^{3}$\newline
\newline
\begin{tabular}{ll}
\begin{tabular}{l}
Newton polygon: \\
\\
\\
\\
\end{tabular}
& \scalebox{.33}{\includegraphics{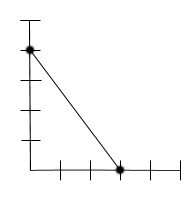}}%
\end{tabular}
\newline
\newline
$%
y^{3}+ax^{4}+bx^{3}y+cx^{2}y^{2}+dxy^{3}+ey^{4}+fx^{5}+gx^{4}y+hx^{3}y^{2}+jx^{2}y^{3}+kxy^{4}+ly^{5}=0,a\neq
0 $.
\newline
Puiseux jets: \newline
$y=(-a)^{1/3}x^{4/3}$.\newline

\begin{tabular}{ll}
\begin{tabular}{l}
10. Diagram: \\
\\
\\
\\
\end{tabular}
& \scalebox{.33}{\includegraphics{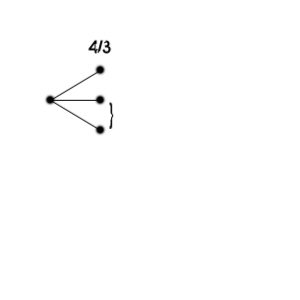}}%
\end{tabular}
\newline

\begin{tabular}{ll}
\begin{tabular}{l}
Newton polygon: \\
\\
\\
\\
\end{tabular}
& \scalebox{.33}{\includegraphics{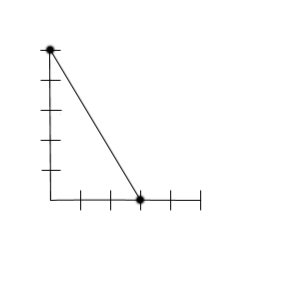}}%
\end{tabular}
\newline
\newline
$%
y^{3}+ax^{5}+bx^{2}y^{2}+cxy^{3}+dy^{4}+ex^{4}y+fx^{3}y^{2}+gx^{2}y^{3}+hxy^{4}+jy^{5}=0,a\neq 0
$.\newline
\newline
Puiseux jets: \newline
$y(-a)^{1/3}x^{5/3}$.\newline
\newline
\begin{tabular}{ll}
\begin{tabular}{l}
11. Diagram: \\
\\
\\
\\
\end{tabular}
& \scalebox{.33}{\includegraphics{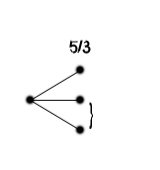}}%
\end{tabular}
\newline
\newline
\begin{tabular}{ll}
\begin{tabular}{l}
Newton polygon: \\
\\
\\
\\
\end{tabular}
& \scalebox{.33}{\includegraphics{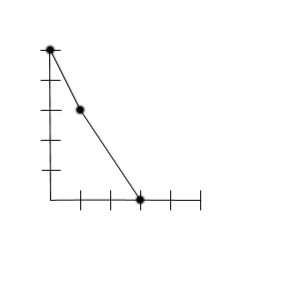}}%
\end{tabular}
\newline
\newline
$%
x^{5}-x^{3}y+y^{3}+ax^{4}y+bx^{2}y^{2}+cx^{3}y^{2}+dxy^{3}+ex^{2}y^{3}+fy^{4}+gxy^{4}+hy^{5}=0
$.\newline
\newline
Puiseux jets: \newline
$y=0x^{3/2}$\newline
$y=x^{3/2}$\newline
\newline
\begin{tabular}{ll}
\begin{tabular}{l}
12. Diagram: \\
\\
\\
\\
\end{tabular}
& \scalebox{.33}{\includegraphics{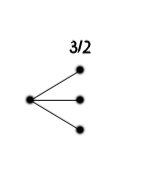}}%
\end{tabular}
\newline
\newline
Tangent cone: $y^{2}(y-x)$.\newline
\newline
\begin{tabular}{ll}
\begin{tabular}{l}
Newton polygon: \\
\\
\\
\\
\end{tabular}
& \scalebox{.33}{\includegraphics{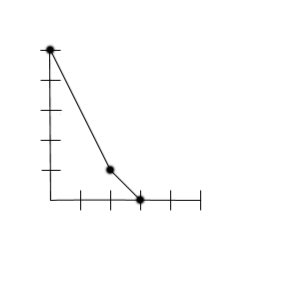}}%
\end{tabular}
\newline
\newline
$%
y^{2}(y-x)+ax^{4}+bx^{3}y+cx^{2}y^{2}+dxy^{3}+ey^{4}+fx^{5}+gx^{4}y+hx^{3}y^{2}+jx^{2}y^{3}+kxy^{4}+ly^{5}=0,a\neq 0
$.\newline
\newline
Puiseux jets: \newline
$y=x$\newline
$y=a^{1/2}x^{3/2}$.\newline
\newline
\begin{tabular}{ll}
\begin{tabular}{l}
13. Diagram: \\
\\
\\
\\
\end{tabular}
& \scalebox{.33}{\includegraphics{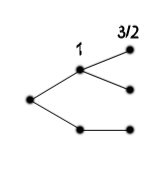}}%
\end{tabular}
\newline
\newline
\begin{tabular}{ll}
\begin{tabular}{l}
Newton polygon: \\
\\
\\
\\
\end{tabular}
& \scalebox{.33}{\includegraphics{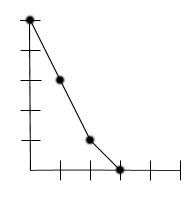}}%
\end{tabular}
\newline
\newline
Quasihomogeneous factors: $x(x^{2}-y)(x^{2}+y)$.\newline
\newline
$%
y^{2}(y-x)+x^{5}+ax^{4}y+bx^{2}y^{2}+cx^{3}y^{2}+dxy^{3}+ex^{2}y^{3}+fy^{4}+gxy^{4}+hy^{5}=0
$.\newline
\newline
Puiseux jets: \newline
$y=x$\newline
$y=x^{2}$\newline
$y=-x^{2}$\newline
\newline
\begin{tabular}{ll}
\begin{tabular}{l}
14. Diagram: \\
\\
\\
\\
\end{tabular}
& \scalebox{.33}{\includegraphics{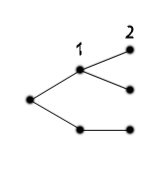}}%
\end{tabular}
\newline
\newline
Quasihomogeneous factors: $x(-y^{2}-x^{4})$.\newline
\newline
$%
y^{2}(y-x)-x^{5}+ax^{4}y+bx^{2}y^{2}+cx^{3}y^{2}+dxy^{3}+ex^{2}y^{3}+fy^{4}+gxy^{4}+hy^{5}=0
$.\newline
\newline
Puiseux jets: \newline
$y=x$\newline
$y=\pm ix^{2}$.\newline
\newline
\begin{tabular}{ll}
\begin{tabular}{l}
15. Diagram: \\
\\
\\
\\
\end{tabular}
& \scalebox{.33}{\includegraphics{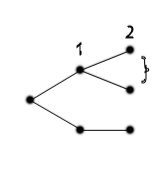}}%
\end{tabular}
\newline
\newline
Quasihomogeneous factors: $-x(y-x^{2})^{2}$.\newline
\newline
$%
y^{2}(y-x)-x^{5}+2x^{3}y+ax^{4}y+bx^{2}y^{2}+cx^{3}y^{2}+dxy^{3}+ex^{2}y^{3}+fy^{4}+gxy^{4}+hy^{5}=0
$.\newline
\newline
Now we use Maple to compute Puiseux expansions.\newline
\newline
Condition: $a+b+1\neq 0$.\newline
\newline
Puiseux jets:\newline
$y=x$\newline
$y=x^{2}+(a+b+1)^{1/2}x^{5/2}$.\newline
\newline
\begin{tabular}{ll}
\begin{tabular}{l}
16. Diagram: \\
\\
\\
\\
\end{tabular}
& \scalebox{.33}{\includegraphics{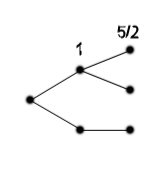}}%
\end{tabular}
\newline
\newline
Case: $b=-a-1$.\newline
\newline
Condition: $d\neq -\frac{1}{4}(1-2a+a^{2}+4c)$.\newline
\newline
Puiseux jets:\newline
$y=x$\newline
$y=x^{2}+x^{3}RootOf(Z^{2}+(a-1)Z-d-c)$.\newline
\newline
\begin{tabular}{lll}
Diagrams: &  &  \\
\begin{tabular}{l}
17. \\
\\
\\
\\
\end{tabular}
& \scalebox{.33}{\includegraphics{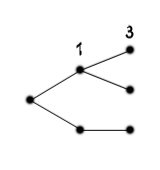}} &
\begin{tabular}{l}
$1-2a+a^{2}+4d+4c>0$ \\
\\
\\
\\
\end{tabular}
\\
\begin{tabular}{l}
18. \\
\\
\\
\\
\end{tabular}
& \scalebox{.33}{\includegraphics{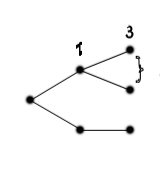}} &
\begin{tabular}{l}
$1-2a+a^{2}+4d+4c<0$ \\
\\
\\
\\
\end{tabular}
\\
&  &
\end{tabular}%
\newline
\newline
Case: $d=-\frac{1}{4}(1-2a+a^{2}+4c)$.\newline
\newline
Condition: $\frac{1}{8}-\frac{1}{8}a-\frac{1}{8}a^{2}-\frac{1}{2}c+\frac{1}{2%
}ac+\frac{1}{8}a^{3}+e+f\neq 0$.\newline
\newline
Puiseux jets: \newline
$y=x$\newline
$y=x^{2}+x^{3}(\frac{1}{2}-\frac{1}{2}a)+(\frac{1}{8}-\frac{1}{8}a-\frac{1}{8%
}a^{2}-\frac{1}{2}c+\frac{1}{2}ac+\frac{1}{8}a^{3}+e+f)^{1/2}x^{7/2}$.%
\newline
\newline
\begin{tabular}{ll}
\begin{tabular}{l}
19. Diagram: \\
\\
\\
\\
\end{tabular}
& \scalebox{.33}{\includegraphics{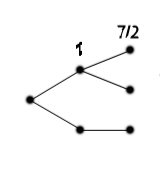}}%
\end{tabular}
\newline
\newline
Case: $f=-(\frac{1}{8}-\frac{1}{8}a-\frac{1}{8}a^{2}-\frac{1}{2}c+\frac{1}{2}%
ac+\frac{1}{8}a^{3}+e)$.\newline
\newline
Condition: $%
D_{1}=256c^{2}-128c+384ca^{2}-256ca+80-32a^{2}-64a+80a^{4}-64a^{3}-512e+1024g+512ae\neq 0
$.\newline
\newline
Puiseux jets: \newline
$y=x$\newline
$y=x^{2}+x^{3}(\frac{1}{2}-\frac{1}{2}%
a)+x^{4}RootOf(16Z^{2}+(16c-20-4a^{2}+24a)Z+8e-14a-16g-8ca^{2}-a^{4}-8c+5+16ac-2a^{3}+12a^{2}-8ea)
$.\newline
\newline
\begin{tabular}{lll}
Diagrams: &  &  \\
\begin{tabular}{l}
20. \\
\\
\\
\\
\end{tabular}
& \scalebox{.33}{\includegraphics{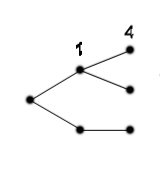}} &
\begin{tabular}{l}
$D_{1}>0$ \\
\\
\\
\\
\end{tabular}
\\
\begin{tabular}{l}
21. \\
\\
\\
\\
\end{tabular}
& \scalebox{.33}{\includegraphics{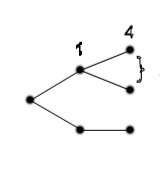}} &
\begin{tabular}{l}
$D_{1}<0$ \\
\\
\\
\\
\end{tabular}%
\end{tabular}%
\newline
\newline
Case: $D_{1}=0$.\newline
\newline
Condition: $D_{2}=\frac{3}{32}-\frac{7a}{64}-\frac{3e}{8}-\frac{c}{8}+h-%
\frac{a^{3}}{32}-\frac{ca}{8}+\frac{ca^{2}}{8}+\frac{ea}{4}+\frac{a^{4}}{32}+%
\frac{a^{3}c}{8}+\frac{a^{2}e}{8}+\frac{c^{2}a}{4}+\frac{ce}{2}+\frac{a^{5}}{%
64}\neq 0$.\newline
\newline
Puiseux jets: \newline
$y=x$ \newline
\newline
$y=x^{2}+(\frac{1}{2}-\frac{1}{2}a)x^{3}+(\frac{a^{2}}{3}+\frac{5}{8}-\frac{c%
}{2}-\frac{3}{4}a)x^{4}+{D_{2}}^{1/2}x^{9/2}$.\newline
\newline
\begin{tabular}{ll}
\begin{tabular}{l}
22. Diagram: \\
\\
\\
\\
\end{tabular}
& \scalebox{.33}{\includegraphics{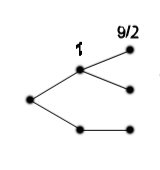}}%
\end{tabular}
\newline

Case: $D_{2}=0$.\newline

The Maple command "factor" applied to the resulting family then
shows that each curve in the family is reducible. Thus, we are
done analyzing the tangent cone $y^{2}(y-x)$.\newline
\newline
Tangent cone: $(y-x)(y-2x)(y-3x)$.\newline
\newline
$%
(y-x)(y-2x)(y-3x)+ax^{4}+bx^{3}y+cx^{2}y^{2}+dxy^{3}+ey^{4}+fx^{5}+gx^{4}y+hx^{3}y^{2}+jx^{2}y^{3}+kxy^{4}+ly^{5}=0
$.\newline

Puiseux jets:\newline
$y=x$\newline
$y=2x$\newline
$y=3x$.\newline
\newline
\begin{tabular}{ll}
\begin{tabular}{l}
23. Diagram: \\
\\
\\
\\
\end{tabular}
& \scalebox{.33}{\includegraphics{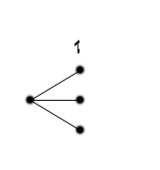}}%
\end{tabular}
\newline
\newline
Tangent cone: $(y-x)(x^{2}+y^{2})$.\newline
\newline
$%
(y-x)(x^{2}+y^{2})+ax^{4}+bx^{3}y+cx^{2}y^{2}+dxy^{3}+ey^{4}+fx^{5}+gx^{4}y+hx^{3}y^{2}+jx^{2}y^{3}+kxy^{4}+ly^{5}=0
$.

Puiseux jets: \newline
$y=x$\newline
$y=\pm ix$.\newline
\newline
\begin{tabular}{ll}
\begin{tabular}{l}
24. Diagram: \\
\\
\\
\\
\end{tabular}
& \scalebox{.33}{\includegraphics{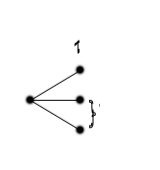}}%
\end{tabular}
\newline
\newline
\underline{Multiplicity 2}\newline
\newline
Tangent cone: $y^{2}$\newline
\newline
\begin{tabular}{ll}
\begin{tabular}{l}
Newton polygon: \\
\\
\\
\\
\end{tabular}
& \scalebox{.33}{\includegraphics{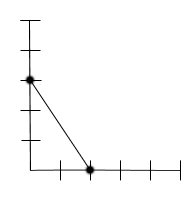}}%
\end{tabular}
\newline
\newline
$%
y^{2}+ax^{3}+bx^{2}y+cxy^{2}+dy^{3}+ex^{4}+fx^{3}y+gx^{2}y^{2}+hxy^{3}+jy^{4}+kx^{5}+lx^{4}y+mx^{3}y^{2}+nx^{2}y^{3}+pxy^{4}+qy^{5}=0,a\neq 0
$.\newline
\newline
Puiseux jets:\newline
$y=(-a)^{1/2}x^{3/2}$.\newline
\newline
\begin{tabular}{ll}
\begin{tabular}{l}
25. Diagram: \\
\\
\\
\\
\end{tabular}
& \scalebox{.33}{\includegraphics{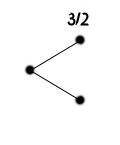}}%
\end{tabular}
\newline
\newline
\begin{tabular}{ll}
\begin{tabular}{l}
Newton polygon: \\
\\
\\
\\
\end{tabular}
& \scalebox{.33}{\includegraphics{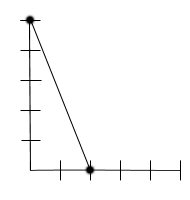}}%
\end{tabular}
\newline
\newline
$%
y^{2}-x^{5}+axy^{2}+by^{3}+cx^{3}y+dx^{2}y^{2}+exy^{3}+fy^{4}+gx^{4}y+hx^{3}y^{2}+jx^{2}y^{3}+kxy^{4}+ly^{5}=0
$.\newline
\newline
Puiseux jets: \newline
$y=\pm x^{5/2}$.\newline
\newline
\begin{tabular}{ll}
\begin{tabular}{l}
26. Diagram: \\
\\
\\
\\
\end{tabular}
& \scalebox{.33}{\includegraphics{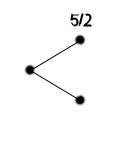}}%
\end{tabular}
\newline
\newline
\begin{tabular}{ll}
\begin{tabular}{l}
Newton polygon: \\
\\
\\
\\
\end{tabular}
& \scalebox{.33}{\includegraphics{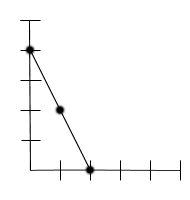}}%
\end{tabular}
\newline
\newline
Quasihomogeneous factors: $(y-x^{2})(y+x^{2})$.\newline
\newline
$%
(y-x^{2})(y+x^{2})+ax^{5}+bx^{3}y+cxy^{2}+dx^{4}y+ex^{2}y^{2}+fy^{3}+gx^{3}y^{2}+hxy^{3}+jx^{2}y^{3}+ky^{4}+lxy^{4}+my^{5}=0
$.\newline
\newline
Puiseux jets:\newline
$y=\pm x^{2}$.\newline
\newline
\begin{tabular}{ll}
\begin{tabular}{l}
27. Diagram: \\
\\
\\
\\
\end{tabular}
& \scalebox{.33}{\includegraphics{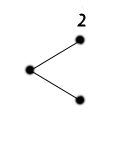}}%
\end{tabular}
\newline
\newline
Quasihomogeneous factors: $y^{2}+x^{4}$.\newline
\newline
$%
y^{2}+x^{4}+ax^{5}+bx^{3}y+cxy^{2}+dx^{4}y+ex^{2}y^{2}+fy^{3}+gx^{3}y^{2}+hxy^{3}+jx^{2}y^{3}+ky^{4}+lxy^{4}+my^{5}=0
$.\newline
\newline
Puiseux jets:\newline
$y=\pm ix^{2}$.\newline
\newline
\begin{tabular}{ll}
\begin{tabular}{l}
28. Diagram: \\
\\
\\
\\
\end{tabular}
& \scalebox{.33}{\includegraphics{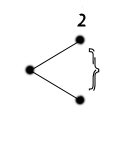}}%
\end{tabular}
\newline
\newline
Quasihomogeneous factors: $(y+x^{2})^{2}$.\newline
\newline
In order to find all remaining singular points corresponding to the tangent
cone $y^{2}$, we should consider the family \newline
\newline
$%
(y+x^{2})^{2}+ax^{5}+bx^{3}y+cxy^{2}+dx^{4}y+ex^{2}y^{2}+fy^{3}+gx^{3}y^{2}+hxy^{3}+jx^{2}y^{3}+ky^{4}+lxy^{4}+my^{5}=0
$.\newline
\newline
We now perform a sequence of Maple calculations of Puiseux expansions for
this family.\newline
\newline
Condition: $b-a-c\neq 0$.\newline
\newline
Puiseux jets: \newline
$y=-x^{2}+(b-a-c)^{1/2}x^{5/2}$.\newline
\newline
\begin{tabular}{ll}
\begin{tabular}{l}
26.Diagram: \\
\\
\\
\\
\end{tabular}
& \scalebox{.33}{\includegraphics{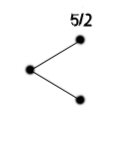}}%
\end{tabular}
\newline
\newline
Case: $b=a+c$.\newline
\newline
Condition: $a^{2}-2ac+c^{2}-4e+4f+4d\neq 0$.\newline
\newline
Puiseux jets: \newline
$y=-x^{2}+x^{3}RootOf(Z^{2}+(a-c)Z+e-f-d)$.\newline
\newline
\begin{tabular}{lll}
Diagrams: &  &  \\
\begin{tabular}{l}
29. \\
\\
\\
\\
\end{tabular}
& \scalebox{.33}{\includegraphics{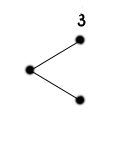}} &
\begin{tabular}{l}
$a^{2}-2ac+c^{2}-4e+4f+4d>0$ \\
\\
\\
\\
\end{tabular}
\\
\begin{tabular}{l}
30. \\
\\
\\
\\
\end{tabular}
& \scalebox{.33}{\includegraphics{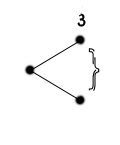}} &
\begin{tabular}{l}
$a^{2}-2ac+c^{2}-4e+4f+4d<0$ \\
\\
\\
\\
\end{tabular}%
\end{tabular}%
\newline
\newline
Case: $e=1/4(a^{2}-2ac+c^{2}+4f+4d)$.\newline
\newline
Condition: $g\neq \frac{a^{2}c}{2}-\frac{ac^{2}}{4}-\frac{a^{3}}{4}+\frac{fa%
}{2}-\frac{fc}{2}-\frac{da}{2}+\frac{dc}{2}+h$.\newline
\newline
Puiseux jets: \newline
$y=-x^{2}+x^{3}(-\frac{a}{2}+\frac{c}{2})+(\frac{a^{2}c}{2}-\frac{ac^{2}}{4}-%
\frac{a^{3}}{4}+\frac{fa}{2}-\frac{da}{2}+\frac{dc}{2}-g+h)^{1/2}x^{7/2}$.%
\newline
\newline
\begin{tabular}{ll}
\begin{tabular}{l}
31. Diagram: \\
\\
\\
\\
\end{tabular}
& \scalebox{.33}{\includegraphics{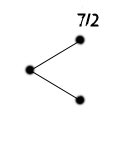}}%
\end{tabular}
\newline
\newline
Case: $g=\frac{a^{2}c}{2}-\frac{ac^{2}}{4}-\frac{a^{3}}{4}+\frac{fa}{2}-%
\frac{fc}{2}-\frac{da}{2}+\frac{dc}{2}+h$.\newline
\newline
Condition: $%
D_{1}=256f^{2}+256fc^{2}-256fa^{2}-512df+256c^{2}a^{2}+256a^{4}+512da^{2}+256d^{2}-512dca+1024j-512hc+512ha-512a^{3}c-1024k\neq 0
$.\newline
\newline
Puiseux jets: \newline
$y=-x^{2}+x^{3}(-\frac{a}{2}+\frac{c}{2}%
)+x^{4}RootOf(16Z^{2}+(16f+8c^{2}-8a^{2}-16d)Z+8dca-16j-4dc^{2}+8hc-6c^{2}a^{2}+16k-8ha-3a^{4}+c^{4}-4da^{2}+8a^{3}c)
$.\newline
\newline
\begin{tabular}{lll}
Diagrams: &  &  \\
\begin{tabular}{l}
32. \\
\\
\\
\\
\end{tabular}
& \scalebox{.33}{\includegraphics{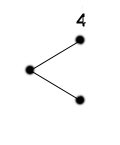}} &
\begin{tabular}{l}
$D_{1}>0$ \\
\\
\\
\\
\end{tabular}
\\
\begin{tabular}{l}
33. \\
\\
\\
\\
\end{tabular}
& \scalebox{.33}{\includegraphics{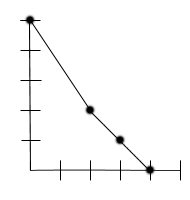}} &
\begin{tabular}{l}
$D_{1}<0$ \\
\\
\\
\\
\end{tabular}%
\end{tabular}%
\newline
\newline
Case: $D_{1}=0$ (solved for $k$).\newline
\newline
Condition: $%
D_{2}=-4adf+8l+4dc^{2}a+4fca^{2}-4dca^{2}-4jc-4hca-fc^{3}+2cf^{2}-2c^{3}a^{2}-2cd^{2}+2af^{2}+2c^{2}a^{3}+2ad^{2}+2hc^{2}-4hf+4hd\neq 0
$.\newline
\newline
Puiseux jets (in parametric form):\newline
$x=-\frac{1}{8}D_{2}T^{2}$\newline
$y=-\frac{1}{64}{D_{2}}^{2}T^{4}-\frac{1}{1024}c{D_{2}}^{3}T^{6}-\frac{1}{%
16384}(c^{2}+2f+2ca-2d){D_{2}}^{4}T^{8}-\frac{1}{32768}{D_{2}}^{5}T^{9}$.%
\newline
\newline
\begin{tabular}{ll}
\begin{tabular}{l}
34. Diagram: \\
\\
\\
\\
\end{tabular}
& \scalebox{.33}{\includegraphics{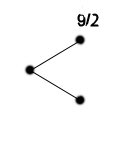}}%
\end{tabular}
\newline
\newline
Case: $D_{2}=0$ (solved for $l$).\newline
\newline
Condition: $%
D_{3}=20480c^{2}a^{2}d+16384m-4096fc^{3}a-8192fc^{2}a^{2}+8192hc^{2}a+4096fc^{2}d-8192hcd-16384d^{2}ca-8192jca+4096a^{2}d^{2}+4096c^{2}a^{4}+4096h^{2}+4096a^{2}f^{2}+8192adh-8192ca^{2}h-8192haf+4096d^{3}-1024f^{2}c^{2}+4096f^{2}d+1024d^{2}c^{2}-8192d^{2}f-8192jf+8192jd+1024c^{4}a^{2}-8192c^{3}a^{3}+16384fcad-2048dc^{3}a-8192a^{3}dc-8192a^{2}df+8192ca^{3}f\neq 0
$.\newline
\newline
Puiseux jets:\newline
$y=-x^{2}+x^{3}(\frac{1}{2}c)+x^{4}(-\frac{1}{4}c^{2}-\frac{1}{2}f-\frac{1}{2%
}ca+\frac{1}{2}%
d)+x^{5}RootOf(64Z^{2}+(64dc+64ad-64ca^{2}-64ac^{2}+64h-64af-128fc-16c^{3})Z+(-144c^{2}a^{2}d-64m+88fc^{3}a+128fc^{2}a^{2}-64hc^{2}a-80fc^{2}d+64hcd-64hcf+96d^{2}ca+32jca+64f^{2}ca-16d^{3}+c^{6}+68f^{2}c^{2}-16f^{2}d+16fc^{4}+12d^{2}c^{2}+32d^{2}f-8hc^{3}+32jf-32jd+20c^{4}a^{2}+64c^{3}a^{3}+8ac^{5}-8dc^{4}-160fcad-32dc^{3}a))
$.\newline
\newline
\begin{tabular}{lll}
Diagrams: &  &  \\
\begin{tabular}{l}
35. \\
\\
\\
\\
\end{tabular}
& \scalebox{.33}{\includegraphics{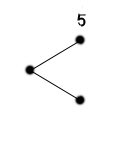}} &
\begin{tabular}{l}
$D_{3}>0$ \\
\\
\\
\\
\end{tabular}
\\
\begin{tabular}{l}
36. \\
\\
\\
\\
\end{tabular}
& \scalebox{.33}{\includegraphics{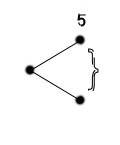}} &
\begin{tabular}{l}
$D_{3}<0$ \\
\\
\\
\\
\end{tabular}%
\end{tabular}%
\newline
\newline
Case: $D_{3}=0$ (solved for $m$).\newline
\newline
Condition: $%
D_{4}=72c^{2}a^{3}d+8afc^{2}d-56ahcd+24ahcf+8a^{3}d^{2}+8c^{2}a^{5}+8ah^{2}+8a^{3}f^{2}+14c^{4}a^{3}-28c^{3}a^{4}-a^{2}c^{5}-12ch^{2}-3f^{2}c^{3}-d^{2}c^{3}+72fca^{2}d-8fc^{3}a^{2}-40fc^{2}a^{3}+40hc^{2}a^{2}-60d^{2}ca^{2}-12f^{2}ca^{2}+16a^{2}dh-16ca^{3}h-16ha^{2}f-36dc^{3}a^{2}-16a^{4}dc-16a^{3}df+16ca^{4}f-4fc^{3}d+8hc^{2}d+12hc^{2}f-14af^{2}c^{2}+4afc^{4}+30ad^{2}c^{2}-8ahc^{3}+2adc^{4}-32ad^{2}f+8dcf^{2}+16daf^{2}-16dhf-8cd^{3}+16ad^{3}+16hd^{3}+16jh-16jaf+16jad-16jca^{2}-8jfc-8jdc+8jac^{2}\neq 0
$.\newline
\newline
Puiseux jets (in parametric form):\newline
$x=-\frac{1}{32}D_{4}T^{2}$\newline
$y=-\frac{1}{1024}{D-2}^{2}T^{4}-\frac{1}{65536}c{D_{4}}^{3}T^{6}-\frac{1}{%
4194304}(c^{2}+2f+2ca-2d){D_{4}}^{4}T^{8}-\frac{1}{268435456}%
(4af-4ad+4ca^{2}+8fc-4dc+ac^{2}+c^{3}-4h){D_{4}}^{5}T^{10}+\frac{1}{%
1073741824}{D_{4}}^{6}T^{11}$.\newline
\newline
\begin{tabular}{ll}
\begin{tabular}{l}
37. Diagram: \\
\\
\\
\\
\end{tabular}
& \scalebox{.33}{\includegraphics{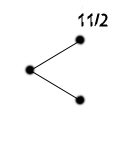}}%
\end{tabular}
\newline
\newline
Case: $D_{4}=0$ (solved for $j$).\newline
\newline
Condition: $%
D_{5}=1024(4a^{2}-28ca+c^{2}+24d-8f)(-2h+2af-2ad+2ca^{2}+fc+dc-ac^{2})^{2}%
\neq 0$.\newline
\newline
(The left side of this condition is the discriminant of the quadratic
polynomial exhibited in the coefficient of $x^{6}$ in the Puiseux jet
immediately below. It is interesting that Maple factored this discriminant.)%
\newline
\newline
Puiseux jets:\newline
$y=-x^{2}+\frac{1}{2}cx^{3}+(-\frac{1}{4}c^{2}-\frac{1}{2}f-\frac{1}{2}ca+%
\frac{1}{2}d)x^{4}+(\frac{1}{2}af-\frac{1}{2}ad+\frac{1}{2}ca^{2}+fc-\frac{1%
}{2}dc+\frac{1}{2}ac^{2}+\frac{1}{8}c^{3}-\frac{1}{2}%
h)x^{5}+x^{6}RootOf(256Z^{2}+(128a^{3}c-160dc^{2}+1024fca-896dca-128a^{2}d-128ah+640c^{2}a^{2}+160ac^{3}+608fc^{2}-448hc-512df+128a^{2}f+32c^{4}256d^{2}+256f^{2})Z+(64d^{4}+c^{8}+10ac^{7}-160a62d^{3}+256fac^{5}-224fdc^{4}+224h^{2}ca-480hc^{3}a^{2}-128hcd^{2}-512hc^{2}a^{3}-256had^{2}-152hac^{4}+144hdc^{3}-576dc^{4}a^{2}+440d^{2}c^{3}a-104dc^{5}a-736a^{4}c^{2}d+608a^{3}cd^{2}-864a^{3}c^{3}d+832d^{2}c^{2}a^{2}-352d^{3}ca-28hc^{5}+192h^{2}c^{2}-96h^{2}d-104d^{3}c^{2}-10dc^{6}+40d^{2}c^{4}+288a^{5}c^{3}+240a^{3}c^{5}+320a^{4}c^{4}+64f^{4}+32fh^{2}+96a^{2}f^{3}-256fd^{3}+38fc^{6}+312f^{3}c^{2}-256f^{3}d+376f^{2}c^{4}+384d^{2}f^{2}+1304f^{2}c^{3}a+1600f^{2}c^{2}a^{2}-696f^{2}c^{2}d-256hcf^{2}+544f^{3}ca-128haf^{2}-352a^{2}df^{2}+480ca^{3}f^{2}+416fa^{2}d^{2}+672fc^{2}a^{4}+424fd^{2}c^{2}-528fhc^{3}+1040fc^{4}a^{2}+1440fc^{3}a^{3}-1568f^{2}cad-2560fc^{2}a^{2}d-1120fhc^{2}a+512fhcd+1376fd^{2}ca+384fadh-640fca^{2}h-1456fdc^{3}a-1088fa^{3}dc+608hdc^{2}a+768hdca^{2}+64c^{6}a^{2}))
$.\newline
\newline
\begin{tabular}{lll}
Diagrams: &  &  \\
\begin{tabular}{l}
38. \\
\\
\\
\\
\end{tabular}
& \scalebox{.33}{\includegraphics{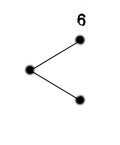}} &
\begin{tabular}{l}
$D_{5}>0$ \\
\\
\\
\\
\end{tabular}
\\
\begin{tabular}{l}
39. \\
\\
\\
\\
\end{tabular}
& \scalebox{.33}{\includegraphics{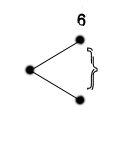}} &
\begin{tabular}{l}
$D_{5}<0$ \\
\\
\\
\\
\end{tabular}%
\end{tabular}%
\newline
\newline
Case: $h=\frac{1}{2}(2af-2ad+2ca^{2}+fc+dc-ac^{2})$.\newline
\newline
The factor command in Maple tells us that each curve in the resulting family
is reducible; in fact, the family becomes\newline
\newline
$\frac{1}{4}(dy+1-acy+ax)(acy^{2}+2x^{2}+cxy+2y+fy^{2}-dy^{2})^{2}$.\newline
\newline
Case: $f=\frac{1}{8}(4a^{2}-28ac+c^{2}+24d)$.\newline
\newline
Condition: $%
D_{6}=512a^{9}+12288c^{3}a^{6}-106176c^{4}a^{5}-122592c^{5}a^{4}-31936c^{6}a^{3}+6144h^{2}a^{3}+768h^{2}c^{3}+32768d^{3}a^{3}+32768d^{3}c^{3}-48c^{6}h+3360c^{7}a^{2}+96c^{7}d-102c^{8}a+24576a^{5}d^{2}+24576c^{2}a^{7}+3072d^{2}c^{5}+340608c^{4}a^{3}d+16896c^{2}a^{5}d-301056c^{2}a^{3}d^{2}-316416c^{3}a^{2}d+97536c^{5}a^{2}d+24576h^{2}ad-27648h^{2}ca^{2}+24576h^{2}dc-26112h^{2}ac^{2}+98304d^{3}ca^{2}+98304d^{3}ac^{2}-36096hc^{2}a^{4}-52034hc^{4}a^{2}+3264hc^{5}a-49152hc^{2}d^{2}-3072hc^{4}d-61440d^{2}ca^{4}+27648ca^{5}h-24576a^{4}dh-49152a^{2}d^{2}h-6432c^{6}ad-49152a^{6}dc+308736c^{3}a^{4}d-4096h^{3}+c^{9}+215040c^{2}a^{2}dh-98304d^{2}c^{4}a+101376hc^{3}ad+86016hcda^{3}-98304hcd^{2}a-118272c^{3}a^{3}h+6144da^{7}-3072a^{6}h-6912a^{8}c\neq 0
$.\newline
\newline
Puiseux jets:\newline
$x=\frac{1}{32768}D_{6}T^{2}$.\newline
$y=-\frac{1}{1073741824}{D_{6}}^{2}T^{4}+\frac{1}{70368744177664}c{D_{6}}%
^{3}T^{6}-\frac{1}{18446744073709551616}(4a^{2}+16d-20ca+5c^{2}){D_{6}}%
^{4}T^{8}+\frac{1}{604462909807314587353088}%
(-8h+4a^{3}+16ad-12ca^{2}+40dc-47ac^{2}+4c^{3}){D-6}^{5}T^{10}-\frac{1}{%
39614081257132168796771975168}%
(8a^{4}+48a^{2}+64d^{2}-44a^{3}c-88dca+34c^{2}a^{2}+112dc^{2}-129ac^{3}+7c^{4}-8ah-28hc)%
{D_{6}}^{6}T^{12}+\frac{1}{40564819207303340847894502572032}{D_{6}}%
^{7}T^{13} $.\newline
\newline
\begin{tabular}{ll}
\begin{tabular}{l}
40. Diagram: \\
\\
\\
\\
\end{tabular}
& \scalebox{.33}{\includegraphics{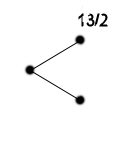}}%
\end{tabular}
\newline
\newline
Case: $D_{6}=0$.\newline
\newline
The factor command in Maple tells us that each curve in the resulting family
is reducible; in fact, the family becomes\newline
\newline
$\frac{1}{256}%
(dy+1-acy+ax)(16x^{2}+16y+4a^{2}y^{2}+16dy^{2}+8cxy-20acy^{2}+c^{2}y^{2})^{2}
$.\newline
\newline
This completes the case of the quasihomogeneous factors $(y+x^{2})^{2}$.%
\newline
\newline
Tangent cone: $y^{2}-x^{2}$.\newline
\newline
\begin{tabular}{ll}
\begin{tabular}{l}
Newton polygon: \\
\\
\\
\\
\end{tabular}
& \scalebox{.33}{\includegraphics{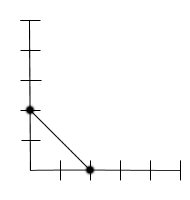}}%
\end{tabular}
\newline
\newline
$y^{2}-x^{2}+$ higher terms $=0$.\newline
\newline
Puiseux jets:\newline
$y=\pm x$.\newline
\newline
\begin{tabular}{ll}
\begin{tabular}{l}
41. Diagram: \\
\\
\\
\\
\end{tabular}
& \scalebox{.33}{\includegraphics{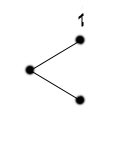}}%
\end{tabular}
\newline
\newline
Tangent cone: $y^{2}+x^{2}$.\newline
\newline
\begin{tabular}{ll}
\begin{tabular}{l}
Newton polygon: \\
\\
\\
\\
\end{tabular}
& \scalebox{.33}{\includegraphics{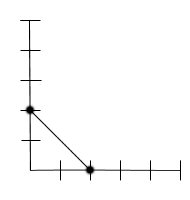}}%
\end{tabular}
\newline
\newline
$y^{2}+x^{2}+$ higher terms $=0$.\newline
\newline
Puiseux jets:\newline
$y=\pm ix$.\newline
\newline
\begin{tabular}{ll}
\begin{tabular}{l}
42. Diagram: \\
\\
\\
\\
\end{tabular}
& \scalebox{.33}{\includegraphics{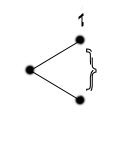}}%
\end{tabular}
\newline
\newline

\section{Reducible curves}

\underline{Factor of degree one}\newline
\newline
If the reducible curve has a factor of degree one, either the line component
passes through the origin or it doesn't. If it does not pass through the
origin, then it is only necessary to list the diagrams from the degree four
case. If the line component does pass through the origin, then, by careful
scrutiny of the Newton polygon, all the corresponding diagrams can easily be
obtained by modifying the diagrams in the preceding case. We now list the
diagrams in this case of a factor of degree one. \newline
\newline
\begin{tabular}{cccccc}
\begin{tabular}{l}
1. \\
\\
\\
\\
\end{tabular}
& \scalebox{.33}{\includegraphics{5fig20}} &
\begin{tabular}{l}
2. \\
\\
\\
\\
\end{tabular}
& \scalebox{.33}{\includegraphics{5fig26}} &
\begin{tabular}{l}
3. \\
\\
\\
\\
\end{tabular}
& \scalebox{.33}{\includegraphics{5fig37}} \\
\begin{tabular}{l}
4. \\
\\
\\
\\
\end{tabular}
& \scalebox{.33}{\includegraphics{5fig38}} &
\begin{tabular}{l}
5. \\
\\
\\
\\
\end{tabular}
& \scalebox{.33}{\includegraphics{5fig60}} &
\begin{tabular}{l}
6. \\
\\
\\
\\
\end{tabular}
& \scalebox{.33}{\includegraphics{5fig62}} \\
\begin{tabular}{l}
7. \\
\\
\\
\\
\end{tabular}
& \scalebox{.33}{\includegraphics{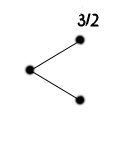}} &
\begin{tabular}{l}
8. \\
\\
\\
\\
\end{tabular}
& \scalebox{.33}{\includegraphics{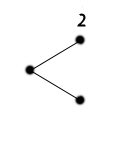}} &
\begin{tabular}{l}
9. \\
\\
\\
\\
\end{tabular}
& \scalebox{.33}{\includegraphics{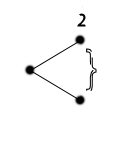}} \\
\begin{tabular}{l}
10. \\
\\
\\
\\
\end{tabular}
& \scalebox{.33}{\includegraphics{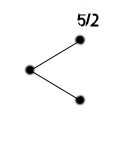}} &
\begin{tabular}{l}
11. \\
\\
\\
\\
\end{tabular}
& \scalebox{.33}{\includegraphics{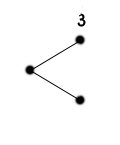}} &
\begin{tabular}{l}
12. \\
\\
\\
\\
\end{tabular}
& \scalebox{.33}{\includegraphics{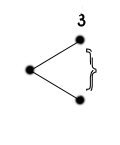}} \\
\begin{tabular}{l}
13. \\
\\
\\
\\
\end{tabular}
& \scalebox{.33}{\includegraphics{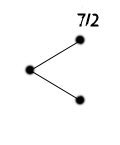}} &
\begin{tabular}{l}
14. \\
\\
\\
\\
\end{tabular}
& \scalebox{.33}{\includegraphics{5fig25}} &
\begin{tabular}{l}
15. \\
\\
\\
\\
\end{tabular}
& \scalebox{.33}{\includegraphics{5fig28}} \\
\begin{tabular}{l}
16. \\
\\
\\
\\
\end{tabular}
& \scalebox{.33}{\includegraphics{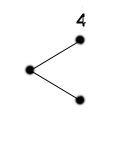}} &
\begin{tabular}{l}
17. \\
\\
\\
\\
\end{tabular}
& \scalebox{.33}{\includegraphics{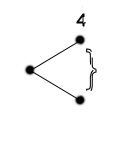}} &
\begin{tabular}{l}
18. \\
\\
\\
\\
\end{tabular}
& \scalebox{.33}{\includegraphics{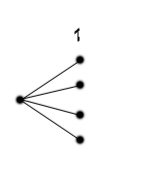}} \\
\begin{tabular}{l}
19. \\
\\
\\
\\
\end{tabular}
& \scalebox{.33}{\includegraphics{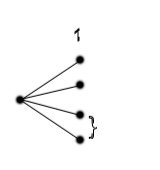}} &
\begin{tabular}{l}
20. \\
\\
\\
\\
\end{tabular}
& \scalebox{.33}{\includegraphics{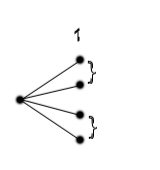}} &  &  \\
&  &  &  &  &
\end{tabular}
\newpage
\begin{tabular}{cccccc}
\begin{tabular}{l}
21. \\
\\
\\
\\
\end{tabular}
& \scalebox{.33}{\includegraphics{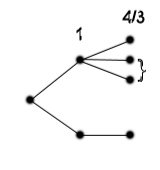}} &
\begin{tabular}{l}
22. \\
\\
\\
\\
\end{tabular}
& \scalebox{.33}{\includegraphics{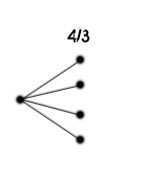}} &
\begin{tabular}{l}
23. \\
\\
\\
\\
\end{tabular}
& \scalebox{.33}{\includegraphics{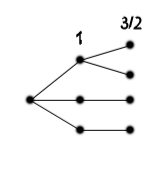}} \\
\begin{tabular}{l}
24. \\
\\
\\
\\
\end{tabular}
& \scalebox{.33}{\includegraphics{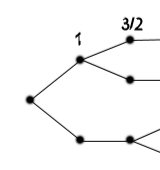}} &
\begin{tabular}{l}
25. \\
\\
\\
\\
\end{tabular}
& \scalebox{.33}{\includegraphics{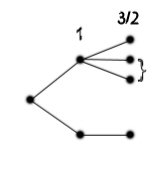}} &
\begin{tabular}{l}
26. \\
\\
\\
\\
\end{tabular}
& \scalebox{.33}{\includegraphics{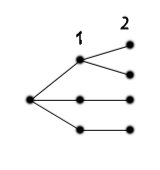}} \\
\begin{tabular}{l}
27. \\
\\
\\
\\
\end{tabular}
& \scalebox{.33}{\includegraphics{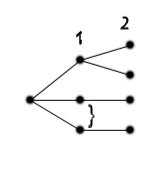}} &
\begin{tabular}{l}
28. \\
\\
\\
\\
\end{tabular}
& \scalebox{.33}{\includegraphics{5fig31}} &
\begin{tabular}{l}
29. \\
\\
\\
\\
\end{tabular}
& \scalebox{.33}{\includegraphics{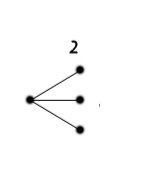}} \\
\begin{tabular}{l}
30. \\
\\
\\
\\
\end{tabular}
& \scalebox{.33}{\includegraphics{5fig29}} &
\begin{tabular}{l}
31. \\
\\
\\
\\
\end{tabular}
& \scalebox{.33}{\includegraphics{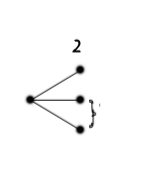}} &
\begin{tabular}{l}
32. \\
\\
\\
\\
\end{tabular}
& \scalebox{.33}{\includegraphics{5fig30}} \\
\begin{tabular}{l}
33. \\
\\
\\
\\
\end{tabular}
& \scalebox{.33}{\includegraphics{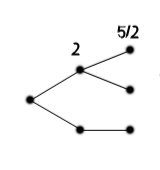}} &
\begin{tabular}{l}
34. \\
\\
\\
\\
\end{tabular}
& \scalebox{.33}{\includegraphics{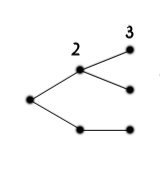}} &
\begin{tabular}{l}
35. \\
\\
\\
\\
\end{tabular}
& \scalebox{.33}{\includegraphics{5fig32}} \\
\begin{tabular}{l}
36. \\
\\
\\
\\
\end{tabular}
& \scalebox{.33}{\includegraphics{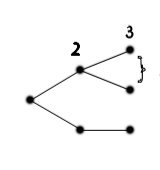}} &
\begin{tabular}{l}
37. \\
\\
\\
\\
\end{tabular}
& \scalebox{.33}{\includegraphics{5fig33}} &
\begin{tabular}{l}
38. \\
\\
\\
\\
\end{tabular}
& \scalebox{.33}{\includegraphics{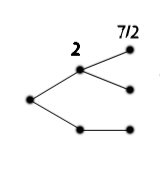}} \\
\begin{tabular}{l}
39. \\
\\
\\
\\
\end{tabular}
& \scalebox{.33}{\includegraphics{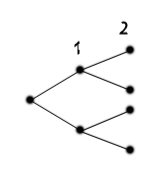}} &
\begin{tabular}{l}
40. \\
\\
\\
\\
\end{tabular}
& \scalebox{.33}{\includegraphics{5fig34}} &
\begin{tabular}{l}
41. \\
\\
\\
\\
\end{tabular}
& \scalebox{.33}{\includegraphics{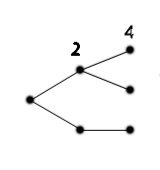}} \\
&  &  &  &  &
\end{tabular}
\newpage
\begin{tabular}{cccccc}
\begin{tabular}{l}
42. \\
\\
\\
\\
\end{tabular}
& \scalebox{.33}{\includegraphics{5fig35}} &
\begin{tabular}{l}
43. \\
\\
\\
\\
\end{tabular}
& \scalebox{.33}{\includegraphics{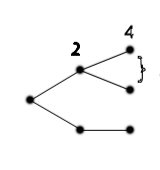}} &
\begin{tabular}{l}
44. \\
\\
\\
\\
\end{tabular}
& \scalebox{.33}{\includegraphics{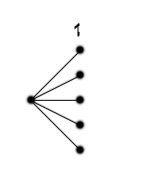}} \\
\begin{tabular}{l}
45. \\
\\
\\
\\
\end{tabular}
& \scalebox{.33}{\includegraphics{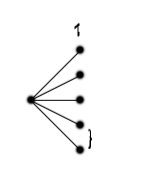}} &
\begin{tabular}{l}
46. \\
\\
\\
\\
\end{tabular}
& \scalebox{.33}{\includegraphics{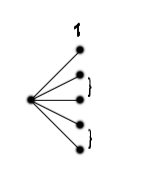}} &  &
\end{tabular}
\newline
\newline
\underline{Factor of degree 2}\newline
\newline
We may assume that the factors of degrees two and three are each
irreducible; otherwise there is a factor of degree one, and we are in the
previous case. If the factor of degree three has a cusp or acnode at the
origin, we obtain nothing new. Furthermore, we only need consider the case
where the tangent line to the conic factor agrees with a tangent line to the
cubic factor.

Some further argument is needed to verify the existence of "complex
conjugate - type" singular points (singular points whose diagrams contain
braces). The Maple computations below do not reveal whether these types of
singular points exist. Real and imaginary parts are inserted for each
coefficient of each factor and then Groebner basis techniques are used to
verify that if the imaginary parts of the coefficients of the product of the
factors are set equal to zero, then there is a solution. Maple is then used
to do a Puiseux expansion of the resulting family, and it is observed that
the singular points in question do or do not exist. There are interesting
phenomena to observe here. Not every diagram without braces is accompanied
by the corresponding diagram with braces. The details of the Groebner basis
computations via Maple are omitted because of their length.[11]\newline
\newline
If the cubic has a crunode at the origin, consider the family\newline
\newline
$(y(y-x)+ax^{3}+bx^{2}y+cxy^{2}+dy^{3})(y+ex^{2}+fxy+gy^{2})=0$\newline
\newline
We will now use Maple to calculate Puiseux expansions.\newline
\newline
Condition: $a\neq -e$.\newline
\newline
Puiseux jets:\newline
$y=x$\newline
$y=ax^{2}$\newline
$y=-ex^{2}$\newline
\newline
\begin{tabular}{lll}
\begin{tabular}{l}
15. and 30. Diagrams: \\
\\
\\
\\
\end{tabular}
& \scalebox{.33}{\includegraphics{5fig28}} & \scalebox{.33}{%
\includegraphics{5fig29}}%
\end{tabular}
\newline
\newline
Case: $a=-e$.\newline
\newline
Condition: $f\neq e-b$.\newline
\newline
Puiseux jets: \newline
$y=x$\newline
$y=-ex^{2}+efx^{3}$\newline
$y=-ex^{2}+(e^{2}-be)x^{3}$.\newline
\newline
\begin{tabular}{lll}
\begin{tabular}{l}
28. and 35. Diagrams: \\
\\
\\
\\
\end{tabular}
& \scalebox{.33}{\includegraphics{5fig31}} & \scalebox{.33}{%
\includegraphics{5fig32}}%
\end{tabular}
\newline
\newline
Case: $f=e-b$.\newline
\newline
Condition: $b\neq e-c-g$.\newline
\newline
Puiseux jets: \newline
$y=x$\newline
$y=-ex^{2}+(e^{2}-be)x^{3}+(3be^{2}-eb^{2}+ce^{2}-2e^{3})x^{4}$\newline
$y=-ex^{2}+(e^{2}-be)x^{3}+(-e^{3}-e^{2}g+2be^{2}-eb^{2})x^{4}$\newline
\newline
\begin{tabular}{lll}
\begin{tabular}{l}
40. and 42. Diagrams: \\
\\
\\
\\
\end{tabular}
& \scalebox{.33}{\includegraphics{5fig34}} & \scalebox{.33}{%
\includegraphics{5fig35}}%
\end{tabular}
\newline
\newline
Case: $b=e-c-g$.\newline
\newline
Condition: $g\neq d$.\newline
\newline
Puiseux jets: \newline
$y=x$\newline
$%
y=-ex^{2}+(ce+eg)x^{3}+(-e^{2}g-eg^{2}-2ecg-ec^{2})x^{4}+(3e^{2}cg+3eg^{2}c+3egc^{2}+ec^{3}+3e^{2}g^{2}+eg^{3})x^{5}
$\newline
$%
y=-ex^{2}+(ce+eg)x^{3}+(-e^{2}g-eg^{2}-2ecg-ec^{2})x^{4}+(e^{3}g-e^{3}d+3e^{2}g^{2}+3e^{2}cg+eg^{3}+3eg^{2}c+3egc^{2}+ec^{3})x^{5}
$.\newline
\newline
\begin{tabular}{lll}
\begin{tabular}{l}
47. Diagram: \\
\\
\\
\\
\end{tabular}
& \scalebox{.33}{\includegraphics{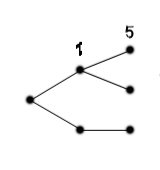}} &
\end{tabular}
\newline
\newline
If $g=d$, then each curve in the resulting family has a linear factor.
\newline
\newline
The final case is where the cubic has a simple point at the origin with
tangent $y=0$.\newline
\newline
$y=(ax^{2}+bxy+cy^{2})(y+dx^{2}+exy+fy^{2}+gx^{3}+hx^{2}y+jxy^{2}+ky^{3})=0$.%
\newline
\newline
Condition: $a\neq d$.\newline
\newline
Puiseux jets:\newline
$y=-dx^{2}$\newline
$y=-ax^{2}$.\newline
\newline
\begin{tabular}{lll}
\begin{tabular}{l}
8. and 9. Diagrams: \\
\\
\\
\\
\end{tabular}
& \scalebox{.33}{\includegraphics{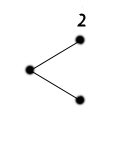}} & \scalebox{.33}{%
\includegraphics{5fig71}}%
\end{tabular}
\newline
\newline
Case: $a=d$\newline
\newline
Condition: $g\neq de-db$.\newline
\newline
Puiseux jets:\newline
$y=-dx^{2}+(de-g)x^{3}$\newline
$y=-dx^{2}+dbx^{3}$.\newline
\newline
\begin{tabular}{lll}
\begin{tabular}{l}
11. and 12. Diagrams: \\
\\
\\
\\
\end{tabular}
& \scalebox{.33}{\includegraphics{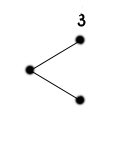}} & \scalebox{.33}{%
\includegraphics{5fig74}}%
\end{tabular}
\newline
\newline
Case: $g=de-db$.\newline
\newline
Condition: $h\neq be+fd-dc-b^{2}$.\newline
Puiseux jets:\newline
$y=-dx^{2}+dbx^{3}+(-d^{2}c-db^{2})x^{4}$\newline
$y=-dx^{2}+dbx^{3}+(-bde-fd^{2}+hd)x^{4}$.\newline
\newline
\begin{tabular}{lll}
\begin{tabular}{l}
16. and 17. Diagrams: \\
\\
\\
\\
\end{tabular}
& \scalebox{.33}{\includegraphics{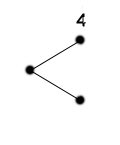}} & \scalebox{.33}{%
\includegraphics{5fig79}}%
\end{tabular}
\newline
\newline
Case: $h=be+fd-dc-b^{2}$.\newline
\newline
Condition: $j\neq bf+ce-2bc$.\newline
\newline
Puiseux jets:\newline
$%
y=-dx^{2}+dbx^{3}+(-d^{2}c-db^{2})x^{4}+(db^{3}+bfd^{2}+d^{2}bc+d^{2}ce-d^{2}j)x^{5}
$.\newline
$y=-dx^{2}+dbx^{3}+(-d^{2}c-db^{2})x^{4}+(3d^{2}bc+db^{3})x^{5}$.\newline
\newline
\begin{tabular}{lll}
\begin{tabular}{l}
48. Diagram: \\
\\
\\
\\
\end{tabular}
& \scalebox{.33}{\includegraphics{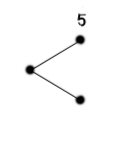}} &
\end{tabular}
\newline
\newline
Case: $j=bf+ce-2bc$.\newline
\newline
Condition: $k\neq cf-c^{2}$.\newline
\newline
Puiseux jets:\newline
$%
y=-dx^{2}+dbx^{3}+(-d^{2}c-db^{2})x^{4}+(3d^{2}bc+db^{3})x^{5}+(-6d^{2}b^{2}c-db^{4}-2d^{3}c^{2})x^{6}
$.\newline
$%
y=-dx^{2}+dbx^{3}+(-d^{2}c-db^{2})x^{4}+(3d^{2}bc+db^{3})x^{5}+(-db^{4}-6d^{2}b^{2}c-d^{3}cf-d^{3}c^{2}+d^{3}k)x^{6}
$.\newline
\newline
\begin{tabular}{lll}
\begin{tabular}{l}
49. Diagram: \\
\\
\\
\\
\end{tabular}
& \scalebox{.33}{\includegraphics{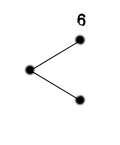}} &
\end{tabular}
\newline
\newline
If $k=cf-c^{2}$, then each curve in the resulting family has a linear factor
and a multiple component. \newline
\newline
This completes the classification.\newline
\newline

\section{Summary of singular points of irreducible real quintic curves}

\begin{tabular}{ccc}
\underline{Multiplicity 2} &  &  \\ \\ \\
\scalebox{.33}{\includegraphics{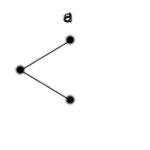}} &
\begin{tabular}{c}
$\frac{2}{2} \leq a \leq \frac{13}{2}$ \\
\\
\\
\\
\end{tabular}
&
\begin{tabular}{c}
$A_1,A_2,A_3,A_4,A_5,A_6,$ \\
$A_7,A_8,A_9,A_{10},A_{11},A_{12}$ \\
\\
\\
\\
\end{tabular}
\\
\scalebox{.33}{\includegraphics{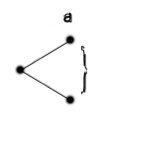}} &
\begin{tabular}{c}
$1 \leq a \leq 6$ \\
\\
\\
\\
\end{tabular}
&
\begin{tabular}{c}
$A^*_1,A^*_3,A^*_5,A^*_7,A^*_9,A^*_{11}$ \\
\\
\\
\\
\end{tabular}
\\
\underline{Multiplicity 3} &  &  \\
\scalebox{.33}{\includegraphics{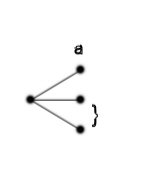}} &
\begin{tabular}{c}
$a=\frac{4}{3},\frac{5}{3},1$ \\
\\
\\
\\
\end{tabular}
&
\begin{tabular}{c}
$E_6,E_8,D^*_4$ \\
\\
\\
\\
\end{tabular}
\\
\scalebox{.33}{\includegraphics{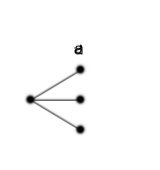}} &
\begin{tabular}{c}
$a=\frac{3}{2},1$ \\
\\
\\
\\
\end{tabular}
&
\begin{tabular}{c}
$E_7,D_4$ \\
\\
\\
\\
\end{tabular}
\\
\scalebox{.33}{\includegraphics{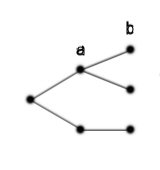}} &
\begin{tabular}{c}
$a=1; b=\frac{3}{2},\frac{4}{2},\frac{5}{2},\frac{6}{2},\frac{7}{2},\frac{8}{%
2},\frac{9}{2}$ \\
\\
\\
\\
\end{tabular}
&
\begin{tabular}{c}
$D_5,D_6,D_7,D_8,D_9,D_{10},D_{11}$ \\
\\
\\
\\
\end{tabular}
\\
&  &
\end{tabular}
\newline
\newpage
\begin{tabular}{ccc}
\scalebox{.33}{\includegraphics{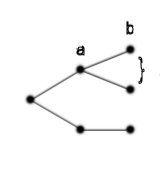}} &
\begin{tabular}{c}
$a=1;b=2,3,4$ \\
\\
\\
\\
\end{tabular}
&
\begin{tabular}{c}
$D^*_6,D^*_8,D^*_{10}$ \\
\\
\\
\\
\end{tabular}
\\
\underline{Multiplicity 4} &  &  \\
\scalebox{.33}{\includegraphics{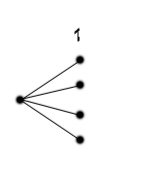}} &  &
\begin{tabular}{c}
$X_9$ \\
\\
\\
\\
\end{tabular}
\\
\scalebox{.33}{\includegraphics{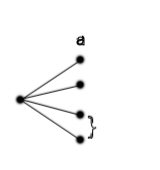}} &
\begin{tabular}{c}
$a=1,\frac{5}{4}$ \\
\\
\\
\\
\end{tabular}
&
\begin{tabular}{c}
$X^*_9,W_{12}$ \\
\\
\\
\\
\end{tabular}
\\
\scalebox{.33}{\includegraphics{5fig82}} &  &
\begin{tabular}{c}
$X^{**}_9$ \\
\\
\\
\\
\end{tabular}
\\
\scalebox{.33}{\includegraphics{5fig4}} &  &
\begin{tabular}{c}
$Z_{11}$ \\
\\
\\
\\
\end{tabular}
\\
\scalebox{.33}{\includegraphics{5fig6}} &  &
\begin{tabular}{c}
$Y^1_{1,1}$ \\
\\
\\
\\
\end{tabular}
\\
&  &
\end{tabular}
\newline
\newpage
\begin{tabular}{ccc}
\scalebox{.33}{\includegraphics{5fig8}} &  &
\begin{tabular}{c}
$Y^{1*}_{1,1}$ \\
\\
\\
\\
\end{tabular}
\\
\scalebox{.33}{\includegraphics{5fig10}} &  &
\begin{tabular}{c}
$X_{1,1}$ \\
\\
\\
\\
\end{tabular}
\\
\scalebox{.33}{\includegraphics{5fig12}} &  &
\begin{tabular}{c}
$X^*_{1,1}$ \\
\\
\\
\\
\end{tabular}
\\
&  &
\end{tabular}%
\newline

\section{Summary of singular points of reducible real quintic curves}

\begin{tabular}{ccc}
\underline{Multiplicity 2} &  &  \\ \\ \\
\scalebox{.33}{\includegraphics{6fig201}} &
\begin{tabular}{c}
$a=\frac{2}{2},\frac{3}{2},\frac{4}{2},\frac{5}{2},\frac{6}{2},\frac{7}{2},%
\frac{8}{2},5,6$ \\
\\
\\
\\
\end{tabular}
&
\begin{tabular}{c}
$A_1,A_2,A_3,A_4,A_5,A_6,$ \\
$A_7,A_9,A_{11}$ \\
\\
\\
\\
\end{tabular}
\\
\scalebox{.33}{\includegraphics{6fig202}} &
\begin{tabular}{c}
$a=1,2,3,4$ \\
\\
\\
\\
\end{tabular}
&
\begin{tabular}{c}
$A^*_1,A^*_3,A^*_5,A^*_7$ \\
\\
\\
\\
\end{tabular}
\\
\underline{Multiplicity 3} &  &  \\
\scalebox{.33}{\includegraphics{6fig205}} &
\begin{tabular}{c}
$a=1,\frac{3}{2},2$ \\
\\
\\
\\
\end{tabular}
&
\begin{tabular}{c}
$D_4,E_7,J_{10}$ \\
\\
\\
\\
\end{tabular}
\\
\scalebox{.33}{\includegraphics{6fig206}} &
\begin{tabular}{c}
$a=1,\frac{4}{3},2$ \\
\\
\\
\\
\end{tabular}
&
\begin{tabular}{c}
$D^*_4,E_6,J^*_{10}$ \\
\\
\\
\\
\end{tabular}
\\
\scalebox{.33}{\includegraphics{6fig203}} &
\begin{tabular}{c}
$a=1;b=\frac{3}{2},\frac{4}{2},\frac{5}{2},\frac{6}{2},\frac{7}{2},\frac{8}{2%
},5$ or, \\
\\
$a=2;b=\frac{5}{2},\frac{6}{2},\frac{7}{2},\frac{8}{2}$ \\
\\
\\
\\
\end{tabular}
&
\begin{tabular}{c}
$D_5,D_6,D_7,D_8,D_9,D_{10},D_{12}$ \\
\\
$J_{11},J_{12},J_{13},J_{14}$ \\
\\
\\
\\
\end{tabular}
\\
&  &
\end{tabular}%
\newline
\newpage
\begin{tabular}{ccc}
\scalebox{.33}{\includegraphics{6fig204}} &
\begin{tabular}{c}
$a=1;b=2,3,4$ or, \\
\\
$a=2;b=3,4$ \\
\\
\\
\\
\end{tabular}
&
\begin{tabular}{c}
$D^*_6,D^*_8,D^*_{10}$ \\
\\
$J^*_{12},J^*_{14}$ \\
\\
\\
\\
\end{tabular}
\\
\underline{Multiplicity 4} &  &  \\
\scalebox{.33}{\includegraphics{5fig89}} &  &
\begin{tabular}{c}
$X_9$ \\
\\
\\
\\
\end{tabular}
\\
&  &  \\
\scalebox{.33}{\includegraphics{6fig222}} &
\begin{tabular}{c}
$a=1,\frac{4}{3}$ \\
\\
\\
\\
\end{tabular}
&
\begin{tabular}{c}
$X^*_9,W_{13}$ \\
\\
\\
\\
\end{tabular}
\\
\scalebox{.33}{\includegraphics{5fig82}} &  &
\begin{tabular}{c}
$X^{**}_9$ \\
\\
\\
\\
\end{tabular}
\\
\scalebox{.33}{\includegraphics{5fig4}} &  &
\begin{tabular}{c}
$Z_{11}$ \\
\\
\\
\\
\end{tabular}
\\
\scalebox{.33}{\includegraphics{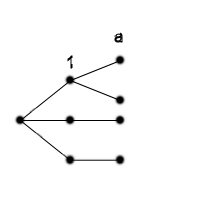}} &
\begin{tabular}{c}
$a=\frac{3}{2},2$ \\
\\
\\
\\
\end{tabular}
&
\begin{tabular}{c}
$X_{1,1},X_{1,2}$ \\
\\
\\
\\
\end{tabular}
\\
&  &
\end{tabular}
\newline
\newpage
\begin{tabular}{ccc}
\scalebox{.33}{\includegraphics{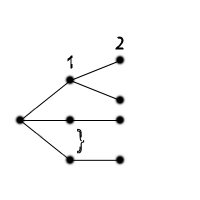}} &  &
\begin{tabular}{c}
$X^{*1}_{1,2}$ \\
\\
\\
\\
\end{tabular}
\\
\scalebox{.33}{\includegraphics{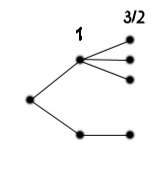}} &  &
\begin{tabular}{c}
$Z_{12}$ \\
\\
\\
\\
\end{tabular}
\\
\scalebox{.33}{\includegraphics{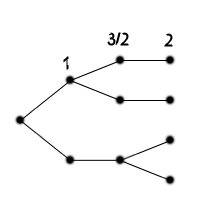}} &  &
\begin{tabular}{c}
$Y^1_{1,2}$ \\
\\
\\
\\
\end{tabular}
\\
\scalebox{.33}{\includegraphics{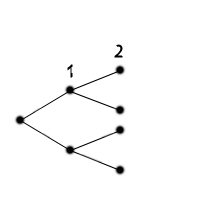}} &  &
\begin{tabular}{c}
$Y^1_{2,2}$ \\
\\
\\
\\
\end{tabular}
\\
\underline{Multiplicity 5} &  &  \\
\scalebox{.33}{\includegraphics{5fig118}} &
\begin{tabular}{c}
$N_{16}$ \\
\\
\\
\\
\end{tabular}
&  \\
\scalebox{.33}{\includegraphics{5fig119}} &
\begin{tabular}{c}
$N^*_{16}$ \\
\\
\\
\\
\end{tabular}
&  \\
\scalebox{.33}{\includegraphics{5fig120}} &
\begin{tabular}{c}
$N^{**}_{16}$ \\
\\
\\
\\
\end{tabular}
&  \\
&  &
\end{tabular}%
\newline

\renewcommand{\baselinestretch}{1}

\textsc{David A. Weinberg: Department of Mathematics and Statistics, Texas
Tech University, Lubbock, TX 79409-1042}

e-mail address: david.weinberg@ttu.edu

\bigskip

\textsc{Nicholas J. Willis: Department of Mathematics and Computer Science,
Whitworth University, Spokane, WA 99251}

e-mail address: nwillis@whitworth.edu\textsc{\ }

\end{document}